\documentclass[12pt]{amsart}
\usepackage{euscript}           

\usepackage{amsmath,amsthm}     
\usepackage{amssymb}            
\usepackage{mathrsfs}           
\usepackage[usenames,dvipsnames]{color}
\usepackage{graphicx}
\usepackage{tikz-cd}            
\usetikzlibrary{shapes}

\usepackage{comment} 
\usepackage{stmaryrd}
\usepackage{hyperref}
\usepackage{cleveref}
\usepackage{manyfoot} 

\hypersetup{
    colorlinks=true,
    linkcolor=blue,
    filecolor=blue, 
    urlcolor=cyan,
}

\newlength{\labwidth}
\newcommand{\labarrow}[1]{
\settowidth{\labwidth}{$\scriptstyle \;\; #1 \;\;$}
\stackrel{#1}{\smash{\hbox to \labwidth{\rightarrowfill}}
\vphantom{\longrightarrow}}
}

\newcommand{\A}{{\mathbb A}}

\newcommand{\C}{{\mathbb C}}

\newcommand{\E}{{\mathbb E}}

\newcommand{\G}{{\mathbb G}}

\renewcommand{\P}{{\mathbb P}}

\newcommand{\R}{{\mathbb R}}
\renewcommand{\S}{{\mathbb S}}

\newcommand{\Z}{{\mathbb Z}}

\newcommand{\calA}{{\mathcal{A}}}

\newcommand{\calC}{{\mathcal{C}}}

\newcommand{\calE}{{\mathcal{E}}}

\newcommand{\calL}{{\mathcal{L}}}
\newcommand{\calM}{{\mathcal{M}}}
\newcommand{\calN}{{\mathcal{N}}}
\newcommand{\calO}{{\mathcal{O}}}

\newcommand{\calS}{{\mathcal{S}}}

\newcommand{\llb}{\llbracket}
\newcommand{\rrb}{\rrbracket}
\newcommand{\ppb}{(\!(q)\!)}
\newcommand{\bbq}{\llbracket q \rrbracket }
\newcommand{\Mellb}{\overline{\mathcal{M}}_{ell}}
\newcommand{\Mell}{\mathcal{M}_{ell}}

\DeclareMathOperator{\TV}{TV}
\DeclareMathOperator{\CAlg}{\mathrm{CAlg}}
\DeclareMathOperator{\SDM}{\mathrm{SDM}}
\DeclareMathOperator{\Fun}{\mathrm{Fun}}
\DeclareMathOperator{\colim}{colim}
\DeclareMathOperator{\Spec}{Spec}

\newtheorem{thm}{Theorem}[section]
\newcounter{numerierer}
\newcounter{leer}

\newtheoremstyle{mydefinition}%
{2ex}{2ex}{\normalfont}{}{\bfseries}{.}{-1.75em}{}

\newtheorem{defn}[thm]{Definition}
\newtheorem{prop}[thm]{Proposition}

\newtheorem{cor}[thm]{Corollary}
\newtheorem{lemma}[thm]{Lemma}

\theoremstyle{definition}  
\newenvironment{definition}{\begin{defn}\rm}{\end{defn}}

\newtheorem{example}[thm]{Example}

\newtheorem{examples}[thm]{Examples}
\newtheorem{construction}{Construction}

\newtheorem{remark}[thm]{Remark}
\newtheorem{warning}[thm]{Warning}
\newtheorem*{convention}{Convention}
\newtheorem*{Question}{Question}

\theoremstyle{mydefinition}
\newtheorem{numex}[thm]{Examples}
\newtheorem{numrmk}[thm]{Remark}

\usepackage[frame,matrix,curve, arrow]{xy}

\xymatrixcolsep{1.9pc}                          
\xymatrixrowsep{1.9pc}
\newdir{ >}{{}*!/-5pt/\dir{>}}                  
\newdir{ |}{{}*!/-5pt/\dir{|}}                  

\subjclass{}

\newcommand{\Sig}{\Sigma_{+}^{\infty}}
\DeclareMathOperator*{\into}{\hookrightarrow}
\DeclareMathOperator*{\Spf}{Spf}

\begin{document}

\title[Tate curve]{A toric construction of the spectral Tate curve}
\author{Laurent Smits}

\date{\today}

\begin{abstract}
    We develop toric geometry spectral algebraic geometry. In particular, we lift toric varieties to flat schemes over the sphere spectrum.
    Furthermore, we prove several rigidity and strictness results allowing us to define and produce strict group actions of flat algebraic tori.
    As an application, we give a construction of the Tate curve in spectral algebraic geometry.
    This approach lets us produce a lift of the Tate curve to $\S\llb q\rrb$    . 
    We furthermore describe how to lift (arithmetic) toric stacks.  
\end{abstract}

\maketitle

\tableofcontents

\section{Introduction and statement of results}

Toric varieties are some of the best understood objects in algebraic geometry. The (purely combinatorial) classification and the structure obtained by the torus action make them a nice family of examples to test theorems on. 
In homotopy theory, one would like to lift commutative rings ``over the sphere'', i.e. for a commutative ring $R$, one would like to have $\E_{\infty}$-rings $\S_R,$ such that $\S_R\otimes_{\S} \Z\cong R$. Rings that admit such a lift can be endowed with particularly nice structure, see \cite[Section IV.1]{NS}.
One can ask a similar question for schemes, namely:\\
\textbf{Question:} Which (flat) schemes over $\Z$ can be obtained as pullback of (flat) Schemes over $\S$ along $\S\to \Z$? 
In this paper, we show that toric schemes do admit such lifts.
\begin{thm}
    Let $A$ be a fixed $\E_{\infty}$-ring. 
    There exists a functor $$\TV^A:\rm{Fans} \to SpSch^\flat_{/\Spec A},$$ into spectral schemes that are flat over $\Spec A$ where 
    \begin{itemize}
        \item $\rm{Fans}$ is a 1-category,
        \item the case $A=\C$ recovers classical toric geometry,
        \item  and the functors $\TV^-$ are compatible with base change. 
    \end{itemize}  
\end{thm}

This theorem is another example of the very rigid structure that toric varieties have. 
Furthermore, we can observe that the functors $\TV^{-}$ commute with certain limits, in particular products, so constructing and working with toric varieties even in spectral algebraic geometry often reduces to 1-categories. 
In \Cref{stacks} and \Cref{arithmetic} we describe extensions of this functor to encompass toric stacks and arithmetic toric varieties, giving, for example, lifts to weighted projective spaces.
We do this by encorporating an auxillary lattice, respectively a group action on the fan, to keep track of either a subgroup of the torus or a twisting group action on the toric variety. \\
As an application of toric geometry, in \Cref{Tate} we will construct a spectral lift of the Tate curve and its regular locus as suggested in \cite[Section 4.3]{Survey}.

\begin{thm}\label{Thmb}
    Let $T\to \Spec \Z\llb q\rrb$ be the Tate curve. 
    There exists a sheaf of $\E_\infty$-rings $\calO^{top}$, such that $(T,\calO^{top})$ is flat over $\Spec \S\llb q\rrb$.   
\end{thm}

The Tate curve plays a fundamental role in the theory of elliptic curves and hence elliptic cohomology, as it gives rise to the compactification of the stack of elliptic curves and the $q$-expansion map of (topological) modular forms. 
As in \cite[Section 4.3]{Survey} this should allow us to construct a compactification of $\Mell^{or}$ of \cite{Ell2}, the moduli stack of oriented elliptic curves. 
This construction is mimicking the classical story and conjecturally coinciding with the compactification constructed by Goerss-Hopkins-Miller \cite{tmf}. \\
However, our methods also allow us  to produce similar compactifications of $\Mell$ and $\Mell^{pre}$, the moduli stack of (strict) elliptic curves and preoriented elliptic curves constructed in \cite{Ell1} and \cite{Ell2}. \\
Different constructions of the oriented Tate curve away from the cusp are also given in \cite{SmoothK} and \cite{Langwen}.

\textbf{Acknowledgments:} I would like to thank Gerd Laures. Furthermore, I would like to thank Jack Davies, Lennart Meier, Björn Schuster, Leonard Tokic, and Langwen Hui for many helpful discussions 
and Jack Davies for such detailed feedback on my drafts. I also want to thank Elmo Vuorenmaa for prompting me to think about toric stacks and arithmetic toric varieties.

\textbf{LLM statememt:} I did not use any LLM for the work presented here. 

\section{Fans and Toric Varieties}

Toric geometry is a part of algebraic geometry that is described purely by combinatorial data, such as cones and fans. A reference to toric geometry in classical algebraic geometry is \cite{Fulton}. 
This textbook covers many more topics than we will need. In this text, we will focus on the basic notions of toric geometry, somewhat streamlined (and restricted to) what we need for applications in elliptic cohomology. \\
For us, a lattice will mean a finitely generated free abelian group. We will denote the 1-category of lattices and group homomorphisms as $Lat$.

\subsection{Cones and affines}

\begin{definition}
    A \emph{cone} $\sigma$ in the lattice $L$ is a finitely generated submonoid $\sigma$ of $L$. 
    A \emph{morphism of cones} $\psi:\sigma\to\tau$ is given by group homomorphism $\psi:L\to N$, with $\psi(\sigma)\subseteq \tau.$  We denote the 1-category given by cones and morphisms of cones by $\mathrm{Cones}$.
\end{definition}

\begin{remark}
    We will often denote cones just by $\sigma$ and keep in mind, that it is equipped with an embedding $\sigma\into L$. 
\end{remark}

\begin{definition}
    A cone $\sigma\subseteq L$ is called \emph{saturated} if for any $v\in L$, such that $rv\in \sigma$ and $r\in \Z_{\leq0}$ already $v\in \sigma$.  
\end{definition}

\begin{remark}\label{doubleconic}
	It is not enough to assume this condition on generators. Let $\sigma \subseteq \Z^2$ be the cone generated by $(1,1)$ and $(1,-1)$. The element $(2,0)$ is in the cone, which is divisible in the lattice but not in the cone. 
\end{remark}

\begin{convention}
    Unless stated otherwise every cone will be saturated. In particular, we always implicitly replace a cone by its saturated closure inside the lattice.  
\end{convention}

\begin{remark}
    The assumption that all cones should be saturated is purely for convenience and could largely be weakened. Considering the saturated case will be enough for the applications we have in mind.
\end{remark}

\begin{numex}\begin{enumerate}
    \item Every lattice $L$ contains trivial cones $\{0\}$ and $L$. The cone $\{0\}$ will be important later.
    \item The submonoid $\Z_{\geq 0}^n$ inside $\Z^n$ is a cone.
    \item The submonoid generated by 2 and 3 in $\Z$ is a cone, but not saturated. This will be a non-example in the remainder of this text. 
    \item Any group-endomorphism $\Z^n\to \Z^n$ defines a distinct endomorphism of cones $\{0\}\subseteq \Z^n$, even though they restrict to the same morphism of monoids $\{0\}\to \{0\}$. 
\end{enumerate} 
\end{numex}

\begin{definition}
    Let $L$ be a lattice with dual lattice $L^*$ and dual pairing $\langle -,- \rangle$. Let $\sigma$ be a cone in $L$. The \emph{dual cone} $\Check{\sigma}$ is given by $\{v\in L^*\mid \langle u,v\rangle \geq 0 \hspace{2pt} \forall u\in \sigma\}$.\\
    The \emph{real cone} $\sigma_\R$ of $\sigma$ is given by $\R_{\geq 0}\cdot\sigma \subseteq L\otimes_{\mathbb{Z}} \R$. Similarly, we will denote the group completion of $\sigma$ by $\sigma_\Z\subseteq L$.
\end{definition}

\begin{definition}
    A cone $\sigma$ with real cone $\sigma_\R$ is called \emph{strongly convex} if no non-trivial real subspace of $L\otimes_\Z \R$ is contained in $\sigma_\R$. 
\end{definition}

\begin{convention}
    In the following, all cones (except their duals) will be strongly convex. This excludes the full trivial cone $L\into L$ for $L\neq\{0\}$ as an example for the remainder of this text. In particular, every cone sattisfies $\sigma=\sigma_\R\cap L\subseteq L\otimes \R$. 
\end{convention}

\begin{lemma}
    The dual cone construction yields a self-duality $$(\check{-}):\rm{Cones}\to \rm{Cones}^{op},$$ which is given on morphisms by the dual of a lattice homomorphism. 
\end{lemma}

\begin{proof}
    This is a straightforward calculation using supporting hyperplanes, see \cite[page 9]{Fulton}. 
\end{proof}

\begin{construction}\label{AffTV}
Let $A$ be a fixed $\E_\infty$-ring. We obtain a functor:
\begin{align}\label{eq}
    \rm{Cones}^{ssc} \xrightarrow{(\Check{-})} \rm{Cones}^{op}&\xrightarrow{forget} \rm{CMon}(Sets)^{op} \nonumber\\ &\xrightarrow[]{\textit{A}\otimes \Sig} \CAlg_\textit{A}^{op}\xrightarrow[]{\Spec} SDM_{/\Spec \textit{A}}, 
\end{align}
Here, the second arrow is the forgetful functor, only remembering the monoid structure of a cone. The symbol $\otimes$ is understood to be the tensor product over the sphere spectrum, also known as the smash product of spectra. \\
\end{construction}

\begin{remark}
    The functor is well defined even on the whole category $\rm{Cones}$, but the varieties in the image are not as well-behaved if we do not restrict our domain to $\rm{Cones}^{ssc}$. 
\end{remark}

\subsection{Fans and globalization}

The spectral Deligne-Mumford stacks in the essential image of the functor of \Cref{AffTV} will be examples of (affine) toric varieties. 
For our purposes, it will be important to construct non-affine toric varieties. 
These will be built by gluing along toric substacks. 
To give a precise definition we have to introduce the necessary combinatorial framework.

\begin{numex}
    \begin{enumerate}
        \item Let $\sigma=\{0\}$ in $L=\Z^n$. The dual cone $\check{\sigma}$ is given by $L^\vee\cong\Z^n$, the dual lattice. The monoid ring $A[\check{\sigma}]\simeq A[x_1^\pm,...,x_n^\pm]$ can be identified with the flat ring of Laurent polynomials over $A$.
        \item Let $\sigma=\Z^n_{\geq0}$ be a cone in $\Z^n$. After identifying $(\Z^n)^\vee\cong \Z^n$, we have $\check{\sigma}=\Z^n_{\geq0}$. So we obtain an identification $\A^n \simeq \Spec A[\check{\sigma}]$.
        \item Let $\sigma$ be the (saturated) cone of Remark \ref{doubleconic}. The associated stack is the double cone $\Spec \S [x,y,z]/(y^2=xz)$. 
    \end{enumerate}
\end{numex}

The stacks in the image of this functor, by definition, will all be affine. To construct more general stacks, we want to introduce a gluing process. We will introduce this via combinatorial notions using fans. 

\begin{definition}
    Let $\sigma$ be a cone in $L$. We say $\tau\subseteq \sigma$ is a \emph{face} of $\sigma$ if $\tau$ is saturated and for all $v,w\in\sigma$ the sum $v+w$ is an element of $\tau$ if and only if $v,w$ are elements of $\tau$.
\end{definition}

\begin{remark}
    Geometrically, we can say that $\tau$ is a face of $\sigma$ if it is a (saturated) cone and either $\tau=\sigma$ or $\tau_{\R}$ is contained in the boundary of $\sigma_{\R}$ (understood as subset of the vector space spanned by $\sigma_\R$).  
\end{remark}

\begin{numex}\label{faces}
\begin{enumerate}
    \item If $\sigma=\Z_{\geq 0}$ in $L=\Z$, then the only faces of $\sigma$ are $\{0\}$ and $\sigma.$ 
    \item Every face of a (saturated, strongly convex) cone is a (saturated, strongly convex) cone itself.
\end{enumerate}
\end{numex}

\begin{lemma}\label{substack}
    Given a cone $\sigma$ and a face $\tau$ of $\sigma$. The functor defined in \Cref{AffTV} sends $\tau\into \sigma$ to an open embedding of spectral Deligne-Mumford stacks.  
\end{lemma}

\begin{proof}
    By construction the algebra $A\otimes \Sig \Check{\sigma}$ is a localization of $A \otimes\Sig \Check{\tau}$ in the sense of \cite[Remark 1.1.4.2]{SAG}, which can be seen as in \cite[page 4]{Fulton}. 
\end{proof}

\begin{examples}
    As we have seen, every cone $\sigma$ in $L$ has $\{0\}\subseteq L$ as a face. In particular, if $L$ has rank $n$, then the Lemma shows that the stack we constructed from $\sigma$ has an open substack isomorphic to $\G_m^n$.
\end{examples}

\begin{definition}
    A \emph{fan} in a lattice $L$ is a set $F$ of (saturated, strongly convex) cones in $L$, such that any face of a cone in $F$ is also contained in $F$ and the intersection of two cones in $F$ is a face of these cones. 
    A \emph{map of fans} $F\to F'$ is a linear map of lattices $\psi: L\to N$ with the property that for any cone $\sigma$ in $F$ there is a cone $\sigma'$ in $F'$ such that $\psi(\sigma)\subseteq \sigma'$. We denote the discrete category of fans and morphisms of fans by $\mathrm{Fans}$.\\
\end{definition}

\begin{numrmk}
    \begin{enumerate}
        \item The category of $\mathrm{Cones}^{ssc}$ embeds fully faithfully into the category of $\mathrm{Fans}$. 
        \item Often, we will define fans by their maximal cones, i.e., all faces are understood.
    \end{enumerate}
\end{numrmk}

\begin{definition}
    A fan that lies in the essential image of $$\mathrm{Cones}^{ssc}\into \mathrm{Fans}$$ is called $\emph{conic}$.
\end{definition}

\begin{examples}
    In the lattice $\Z^n$, the set of all orthants $(\pm\Z_{\geq 0})^n$ and their faces defines a fan.
\end{examples}

\begin{center}
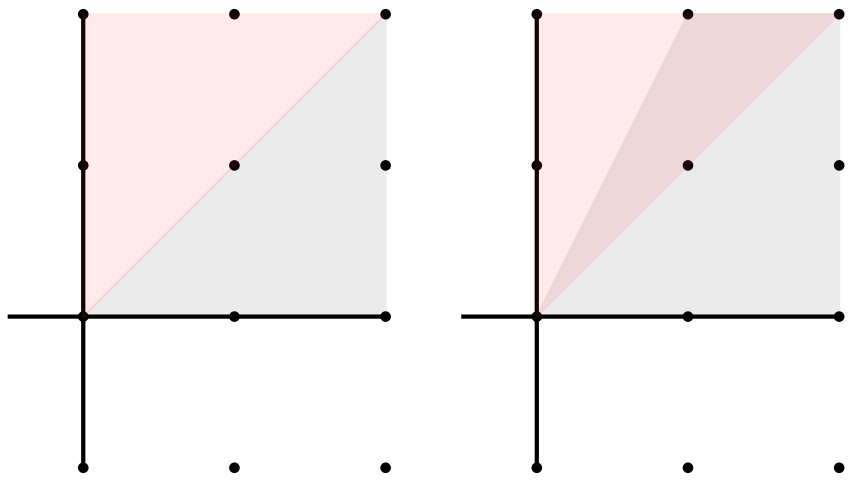
\begin{figure}[h]
\begin{tikzpicture}[scale=0.5]
    \coordinate (a) at (4,0);
    \coordinate (b) at (-1,0);
    \coordinate (c) at (0,4);
    \coordinate (d) at (0,-2);

    \draw [thick] (a) to (b);
    \draw [thick] (c) to (d);

    \fill (2,0) circle (2pt);
    \fill (0,2) circle (2pt);
    \fill (0,4) circle (2pt);
    \fill (2,2) circle (2pt);
    \fill (0,0) circle (2pt);
    \fill (4,0) circle (2pt);
    \fill (4,2) circle (2pt);
    \fill (4,4) circle (2pt);
    \fill (2,4) circle (2pt);
    \fill (0,-2) circle (2pt);
    \fill (2,-2) circle (2pt);
    \fill (4,-2) circle (2pt);

    \filldraw [black!80, opacity=0.1] (0,0)  -- (4,4)  -- (4,0)  -- cycle;
    \filldraw [red!80, opacity=0.1] (0,0)  -- (4,4)  -- (0,4)  -- cycle;

    \coordinate (a) at (10,0);
    \coordinate (b) at (5,0);
    \coordinate (c) at (6,4);
    \coordinate (d) at (6,-2);

    \draw [thick] (a) to (b);
    \draw [thick] (c) to (d);

    \fill (6,0) circle (2pt);
    \fill (6,2) circle (2pt);
    \fill (6,4) circle (2pt);
    \fill (8,2) circle (2pt);
    \fill (8,0) circle (2pt);
    \fill (10,0) circle (2pt);
    \fill (10,2) circle (2pt);
    \fill (10,4) circle (2pt);
    \fill (8,4) circle (2pt);
    \fill (6,-2) circle (2pt);
    \fill (8,-2) circle (2pt);
    \fill (10,-2) circle (2pt);

    \filldraw [black!80, opacity=0.1] (6,0)  -- (8,4)  -- (10,4) -- (10,0)  -- cycle;
    \filldraw [red!80, opacity=0.1] (6,0)  -- (10,4)  -- (6,4)  -- cycle;
\end{tikzpicture}
\caption{An example and a non-example of a fan. The darker area in the second picture is the intersection of the red and the gray cones. This intersection is not a face of either of these cones. }
    \label{fig:enter-label}
\end{figure}
\end{center}

\begin{examples}\label{Projective}
    Let $e_0,...,e_n$ be the standard basis of $\Z^{n+1}$. We denote by $v_0,...,v_n$ the images of $e_0,...,e_n$ under the projection $$\Z^{n+1}\to \Z^{n+1}/(1,...,1)\Z=:L.$$ Any proper subset of $\{v_0,...,v_n\}$ generates a cone in $L$. The set of all these cones is a fan: Let $I, J \subsetneq \{v_0,...,v_n\}$ and $\sigma_I,\sigma_J$ the cones they generate. Every face of the cone $\sigma_I$ can be generated by a subset of $I$. Also we have $\sigma_I\cap\sigma_J =\sigma_{I\cap J}$ for their sections. So both properties of fans are fulfilled. 
\end{examples}

As mentioned before, we want to glue toric stacks. The fans will be our combinatorial description of our ``gluing data''. 
Let us revisit the case of $\P^1$. We can construct $\P^1$ as a pushout of two copies of $\A^1$ identified via inclusions $\G_m\into \A^1$. To obtain a $\P^1$ and not a line with two origins, it is important that exactly one inclusion $\G_m\into \A^1$ is via $z\mapsto z^{-1}$. 
In our combinatorial description, this translates to the following: We want to identify two copies of (the stacks corresponding to) $\Z_{\geq0}\subseteq \Z$ along the inclusions of cones $\{0\}\into \Z_{\geq0}$. The non trivial inclusion is given by the identification $\{0\}\into \Z_{\geq0}\cong \Z_{\leq0}$. 
In particular, we can describe the situation in the following way: We want to consider the pushout of the diagram of cones $\Z_{\leq 0} \xleftarrow{} \{0\} \xrightarrow{} \Z_{\geq0}$ after applying the functor of \Cref{AffTV}. 
More generally, we can observe that $\{\Z_{\leq0},\{0\},\Z_{\geq0}\}$ is a fan, and the arrows in our diagram are the inclusion of faces. In analogy, it is reasonable to ask for general colimits as described by fans. Let us make this precise:\\

By definition, a fan $F$ is a set of cones. We define a partial order on $F$ by considering $\tau\leq \sigma$ if and only if $\tau$ is a face of $\sigma$. 
Given a morphism of fans $\varphi:F\to F'$ we define the map of partially ordered sets $\Tilde{\varphi}:F\to F', \sigma\mapsto \sigma'$, where $\sigma'$ is the minimal cone in $F'$, such that $\varphi(\sigma)\subseteq \sigma'$. \\ 
We obtain a functor $$\iota:\mathrm{Fans}\to \int_{Poset} Fun(-, Cones^{ssc}),$$ by sending each vertex in the partially ordered set $F$ to the cone it is associated to.\\
Recall that for a functor $\mathscr{F}:\calC\to \mathrm{Cat}$ the Grothendieck construction $\int_\calC \mathscr{F}$ consist of objects $(C,X)$, with $C\in \calC$ and $X\in \mathscr{F}(C)$. Morphisms in $\int_\calC \mathscr{F}$ are given as $(f,u):(C,X)\to (D,Y)$ with $(f:C\to D)\in \calC$ and $(u:\mathscr{F}(f)(X)\to Y)\in \mathscr{F}(D)$.\\
In our particular example of interest, the functor $\iota$ sends a fan $F$ to the pair consisting of the poset $F$ as described above and the functor $F\to \mathrm{Cones}^{ssc}$ sending each element of the poset to the cone it is representing.\\
On morphisms the functor $\iota$ sends $\varphi$ to $(\Tilde{\varphi},\varphi_{\mid \sigma})$.\\
Furthermore, we can compose this functor with the functor of \Cref{AffTV} to obtain $$\mathrm{Fans}\to \int_{Poset} \Fun(-, SDM_{/\Spec \textit{A}}).$$
Moreover, the image takes values in the subcategory $\mathcal{L}ax^{colim}(\mathrm{SDM}_{/\Spec \textit{A}})$ of the diagrams in $\mathrm{SDM}_{/\Spec \textit{A}}$ that admit a colimit, as defined by \cite[Notation 3.11]{Colim}. This is because every morphism in a diagram induced by a fan is étale and such diagrams admit colimits, by \cite[Proposition 8.7.6]{Rezksag} which relies on \cite[Warning 1.2.0.3]{SAG} and is proved as \cite[Proposition 21.4.6.4]{SAG}. So we can apply the global colimit functor \cite[Proposition 3.12]{Colim}.

\begin{definition}
    Given an $\E_\infty$-ring $A\in \mathrm{CAlg}$, we define the functor $$\TV_{-}:\mathrm{Fans}\to \SDM_{/ \Spec \textit{A}}$$ as the composition 
    $$\mathrm{Fans}\to \mathcal{L}ax^{colim}(\mathrm{SDM}_{/\Spec \textit{A}})\xrightarrow{colim} \SDM_{/ \Spec \textit{A}},$$ 
    where the first arrow is the functor described above and colim the global colimit functor of \cite[Proposition 3.12]{Colim}. \\
    A spectral Deligne-Mumford stack is called a \emph{toric variety} if it lies in the essential image of $\TV_{-}$. Sometimes we will say that the toric variety is \textit{based on A} if we want to emphasize the $\E_\infty$-ring $A$.
\end{definition}

\begin{warning}
    In the literature, it is common to define toric varieties as schemes with a group action by an algebraic torus with a dense open orbit. As we will see there is a way to reconcile the situation and a natural way to define this action. Furthermore, toric varieties are usually not varieties in the sense of \cite[Definition 1.1.1]{Ell1}.
\end{warning}

\begin{numrmk}\label{schmes}
\begin{enumerate}
    \item $\TV$ refines to a functor over the functor $\Spec$ so we also have functoriality in the chosen $\E_\infty$-ring. See for example the constructions in  \cite{TorSch}.
    \item Affine toric varieties associated to cones will always be schematic. The morphisms in a diagram induced by a fan are open immersions, so the functor $TV_{-}$ factors through spectral schemes. Later we will use the flexibility given by considering stacks, so we do not use this factorization.
\end{enumerate}
\end{numrmk}

\begin{numex}
\begin{enumerate}
    \item If a fan is conic, then the functor $\TV$ agrees with the functor of \Cref{AffTV}, because the fan has a terminal object (as poset) and as such the colimit is described by evaluation at the maximal cone.
    \item Let $F$ be the fan from Example \ref{Projective}. The toric variety of $F$ is, by definition, the (flat) projective space $\P^n_R$ constructed in \cite[Construction 5.4.1.3]{SAG}. 
    \item The toric variety associated to the fan generated by all cones of the form $\pm \Z^n_{\geq0}$ in $\Z^n$ is equivalent to $(\P^1)^n$. 
\end{enumerate}
\end{numex}

\subsection{Properties of $TV$}

\begin{lemma}\label{Allflat}
    Let $A$ be an $\E_\infty$-ring. Every toric variety $\TV_F$ based on $A$ is flat over $\Spec A$.
\end{lemma}

\begin{proof}
    It suffices to check this in the case where the ground ring is the sphere spectrum $\S$. The assertion is local on $\TV_F$. Locally the structure morphism $\TV_F\to \Spec \S$ is given by $\S\to \S[\check{\sigma}]$. This map is flat.
\end{proof}

\begin{cor}
    Let $F$ be a fan. The underlying stack of $\TV_F$ based on $\S$ agrees with the toric variety $\TV_F$ based on $\Z$.
\end{cor}

\begin{proof}
    The cones $\sigma_i$ in $F$ define an \'etale covering of $\TV_{\sigma_i}\to\TV_F$. Both $(-)^\heartsuit$ and $\times \Spec \Z$ preserve coverings. The claim follows, because locally $\S[\check{\sigma}_i]\otimes \Z\cong \pi_0 \S[\check{\sigma}_i]$ yields equivalences $\TV_{\sigma_i}^{\heartsuit}\simeq \TV_{\sigma_i}\times \Spec \Z$. 
\end{proof}

\begin{thm}
    Every toric variety is schematic. 
\end{thm}

\begin{proof}
    Choose an isomorphism $X\simeq \TV_F$ for some fan $F$. Let $\sigma,\sigma'\in F$ be cones, such that the monoid generated by $\sigma$ and $\sigma'$ is another cone $\Sigma$ (which necessarily has $\sigma,\sigma'$ as faces). Let us denote the intersection $\sigma\cap \sigma'$ by $\tau$. We have seen that the maps $R[\check{\Sigma}]\to \R[\check{\sigma}]\to R[\check{\tau}]$ are localizations. Moreover, we can see that \begin{center}
        \begin{tikzcd}
            R[\check{\Sigma}] \ar{r} \ar{d} & R[\check{\sigma}] \ar{d} \\
            R[\check{\sigma}'] \ar{r} & \R[\check{\tau}] 
        \end{tikzcd}
    \end{center} is a pushout in $\E_\infty$-rings. 
    Fix $\sigma\in F$, such that $\sigma$ is not a proper face of any other cone in $F$. We write $X$ for the colimit $\colim_{\tau\in F\setminus \{\sigma\}} \Spec R[\check{\tau}]$ and $U$ for $\colim_{\tau\in F\setminus\{\sigma\}} \Spec R[\check{(\tau\cap \sigma)}]$. By construction of $U,X$ and $\TV_F$, we have a pushout diagram 
    \begin{center}\begin{tikzcd}
        U \ar{r} \ar{d} &X\ar{d} \\ \Spec R[\check{\sigma}] \ar{r} & \TV_F, 
    \end{tikzcd}
    \end{center}
    where $U \to X$ is an open immersion, as explained by the localization description above. Recall that open immersions of spectral Deligne-Mumford stacks are defined by \'etale maps associated to (-1)-truncated objects of the underlying topoi. Proposition 2.5.2.1 (implication (2) to (3)) of \cite{SAG} shows that the pushout of an open immersion is an open immersion, which implies that the map $\Spec R[\check{\sigma}]\to \TV_F$ is an open immersion. Now, because the family $\{\Spec R[\check{\sigma}]\}_{\sigma}$ indexed by the maximal cones covers $\TV_F$, \cite[Proposition 1.6.7.3]{SAG} implies that $\TV_F\simeq X$ is schematic.
\end{proof}

The proof is a generalization, but very analogous to the proof of \cite[Proposition 5.4.1.7]{SAG}, which establishes the claim for the projective space. 
It is a formalization of the intuition we developed before: The (maximal) cones of a fan cover the toric variety, while the smaller cones describe gluings along open immersions.\\

\begin{cor}\label{fansubstack}
    Every toric variety $X\simeq TV_F$ admits a canonical morphism from a torus $\G_m^n$, which is an open immersion, where $n$ is the rank of the lattice $L$ containing $F$.
\end{cor}

\begin{proof}
    This is Lemma \ref{substack} in the case of conic fans. Otherwise, we can follow the previous proof. 
    There, we established already that for maximal cones $\sigma$ the map $\Spec R[\check{\sigma}]\to \TV_F$ is an open immersion. Let $\tau$ be a face of $\sigma$. Then the map $\Spec R[\check{\tau}]\to \Spec R[\check{\sigma}]\to \TV_F$ is an open immersion, because both maps are open immersions. 
\end{proof}

A similar analysis of the coverings lets us compute products of toric varieties.

\begin{lemma}
    The product in the category of fans is given by $$F\times F'=\{\sigma\times \sigma' \subseteq \Z^n \mid \sigma\in F, \sigma'\in F'\}$$ and products commute with $\TV$.
\end{lemma}

\begin{proof}
    The first part is classical and can be easily checked. First assume that $F$ and $F'$ are conic, i.e. given by cones $\sigma$ and $\sigma'$. Then the product is given by $\sigma\times \sigma'$ with dual $\Check{\sigma}\times \Check{\sigma}'$ so the coordinate ring is given by the tensor product and the toric variety as the product. 
    If the fans are not conic, we use the open affine coverings by affine toric varieties induced by the fans. The product $\TV_F\times_{\Spec \textit{A}} \TV_{F'}$ is precisely obtained by gluing together the open substacks $\TV_\sigma \times_{\Spec \textit{A}} \TV_{\sigma'}$ as defined by $F\times F'$, just as in the case of ordinary schemes. The last assertion related to the main input into (the usual constructive proof of) the existence of fiber products of (spectral) Deligne-Mumford stacks, see \cite[Proposition 1.4.11.1]{SAG}.  
\end{proof}

This fact gives us an alternative way to produce the abelian group object structure on $\G_m$. 

\begin{cor}
    The scheme $\G_m=\Spec \S[t^\pm]$ is an abelian group object in $SDM$ via the composition of $TV$ with the functor $Lat^{op}\to \mathrm{Fans}, L\mapsto (\{0\}\into L)$. 
\end{cor}

The category of fans more generally posses all fiber products, but $\TV$ does not commute with all fiber products, which can already be observed based on $\Z$. However, there is a sufficient condition to when $\TV$ commutes with fiberproducts.

\begin{definition}[\cite{Pullback}]
    A morphism of fans $F\to F'$ is called \emph{weakly semistable} if every cone $\sigma\in F$ surjects onto a cone $\tau\in F'$.
\end{definition}

\begin{numex}
    \begin{enumerate}
        \item Projections $F\times F'\to F$ are weakly semistable. 
        \item Let $F_n$ be the fan from Example \ref{Projective} defining the projective $n$-space. The map $F_n\to F_{n-1}$ given by $v_{n-1}\mapsto 0, v_n\mapsto v_{n-1}$ and $v_i\mapsto v_i$ for all $i\in\{0,...,n-2\}$ is weakly semistable.
        \item The map $\Z_{\geq0}\to \Z_{\geq0}$ given by multiplication by $2$ on the lattice $\Z$ is not weakly semistable. 
        \item Let $F$ be the fan in $\Z^2$ defined by the cones $\sigma=\langle (1,0),(1,1)\rangle$ and $\sigma=\langle (0,1),(1,1)\rangle$. Let $F'$ be the conic fan in $\Z^2$ defined by $\Sigma=\langle(1,0),(0,1)\rangle$. The map $F\to F'$ given by the identity on $\Z^2$ is not weakly semistable. 
        \item In \cite{Pullback} it is shown that a morphism is weakly semistable if and only if it induces a flat map with reduced fibers on $\Z$ based toric varieties.
    \end{enumerate}
\end{numex}

\begin{thm}\label{fiberprods}
    The functor $\TV$ preserves fiberproducts $F\times_{F'} G$ if the map $F\to F'$ is weakly semistable. 
\end{thm}

\begin{proof}
    We want to show that the limit comparison map  $$f:\TV_{F \times_{F'}G} \to \TV_F\times_{\TV_{F'}} \TV_G$$ is an equivalence. Using \cite[Corollary 6.1.4.12]{SAG}, this is equivalent to checking that $f$ is locally almost of finite presentation, a relative $n$-stack, and an equivalence after basechanging to some field. \\
    It suffices to check the claim for $A=\S$. Using the spectral Hilbert Basis theorem \cite[Proposition 7.2.4.31]{HA}, we see that $f$ is locally almost of finite presentation. The same holds for $g:\TV_F\to \TV_{F'}$. By \cite[Corollary 6.1.4.9]{SAG}, this implies that $g$ is flat if it is flat after changing the base to a field $k$. 
    This was proven in \cite[Theorem 2.1.4]{Pullback}. Consequently, the stack $\TV_F\times_{\TV_{F'}} \TV_G$ is flat over $\TV_G$ and by Lemma \ref{Allflat} also over $\Spec \S$.\\
    In particular, checking if the map $f:\TV_{F \times_{F'}G} \to \TV_F\times_{\TV_{F'}} \TV_G$ is an equivalence after base change to a field reduces to \cite[Lemma 2.2.6]{Pullback} and is therefore true.\\
    Lastly, $f$ is a relative $0$-stack by \cite[Remark 6.3.1.14]{SAG}, because it is a map of schemes and schemes are closed under products in $SDM$. We therefore checked all conditions of \cite[Corollary 6.1.4.12]{SAG} and conclude that $f$ is an equivalence.
\end{proof}

Let us analyze the structure of toric varieties. In the literature, it is common to define toric varieties as schemes with an open subscheme isomorphic to $\G_m^n$, such that the multiplication of $\G_m^n$ extends to an action on the whole scheme. 
We want to elaborate on what this ``action'' could mean and how it arises.

\begin{definition}
    We define the \emph{action category} as the discrete category $\mathcal{A}ct$ consisting of pairs $(L,b)$ where $L$ is a lattice and $b\in\{0,1\}$ is a boolean and morphisms given by \begin{align*}
        \mathcal{A}ct((A,0),(B,0))=\mathrm{Lat}(A,B)\\
        \mathcal{A}ct((A,1),(B,0))=\mathrm{Lat}(A,B)\\
        \mathcal{A}ct((A,0),(B,1))=\emptyset\\
        \mathcal{A}ct((A,1),(B,1))\subseteq\mathrm{Lat}(A\times\Z,B\times \Z)
    \end{align*}
    where the maps $\mathcal{A}ct((A,1),(B,1))$ consists of the set of lattice morphisms that restrict to the identity on the auxiliary factor $\Z$. Given $(A,1)\xrightarrow{f} (B,1)\xrightarrow{g} (C,0)$ given as lattice morphisms $\varphi$ and $\psi$, then we define the composition $g\circ f$ as lattice morphism $ \psi \circ pr_B \circ \varphi_{\mid A}: A\to C$, where $pr_B:B\times \Z\to B$ is the projection map. In all other cases, the composition is given by composition in $\mathrm{Lat}$.
\end{definition}

\begin{remark}
    The category $Lat^{op}\simeq Lat$ is a full subcategory of $\mathcal{A}ct$, via $L\mapsto (L,0)$.
\end{remark}

\begin{definition}
    We say an abelian group object $A$ acts on an object $X$ if there is a product preserving functor $\Phi:\mathcal{A}ct\to \mathcal{C}$, with $\Phi(0,1)\simeq X$ and that restricts to the abelian group object $A$ via the above embedding. We write $\mathcal{A}ct(\mathcal{C}):=Fun^{\times}(\mathcal{A}ct,\mathcal{C})$ for the total action catgory in $\mathcal{C}$. The morphism $(\Z,1)\to (0,1)$, given by the matrix $(1,1)$ is called action morphism. 
\end{definition}

\begin{remark}\label{actiondiag}
    Such an abelian group action on an object $X$ defines a simplicial diagram called an action diagram. The $n$-th level given by $X\times A^{\times n}$ and the boundary maps are given by the projection onto the object $X$, the action morphism, and multiplication of the abelian group object. The degeneracy maps are given by the unit of $A$. More precisely we obtain a functor $\Delta^{op}\to \calA ct$, which sends the maps $d_0,d_1:[1]\to [0]$ to the matrices $(1\hspace{0.2cm} 1), (0 \hspace{0.2cm} 1)$ understood as morphisms $(\Z,1)\to (\Z,1)$. 
\end{remark}

\begin{lemma}
    Let $F$ be a fan in the lattice $L$. The toric variety $TV_F$ carries a natural torus action, which restricts to the multiplication on the (open) substack given by $\{0\}\into F$ of Corollary \ref{fansubstack}.
\end{lemma}

\begin{proof}
    We have already seen, that $\G_m^n$ embeds into $TV_F$. The map $0\times F\to F$ given by addition on the lattice $L$ defines a product preserving functor $\mathcal{A}ct\to \mathrm{Fans}$, where $\{0\}$ acts on $F$. Postcomposing this with $TV$ yields an action of $\G_m$ on $TV_F$, because $0\times F\to F$ is weakly semistable. Restricting this action to $\G_m^n$ yields, by construction, multiplication.
\end{proof}

\begin{lemma}
    A toric morphism is equivariant with respect to their respective torus actions, in the sense that every toric morphism defines a map in $\mathcal{A}ct(\mathrm{SDM}_{/\Spec A})$.
\end{lemma}

\begin{proof}
    Let $F,F'$ be fans in $L,L'$ and denote the trivial fans in $L,L'$ by $0_L,0_{L'}$. Every map of fans $\varphi:F\to F'$, induces a morphism $0_L\times F\to 0_{L'} \times F'$, which commutes with the action of $0_L$ on $F$ and $0_{L'}$ on $F'$ in $\mathrm{Fans}$. This assembles into a morphism in $\mathcal{A}ct(\mathrm{Fans})$. Again composition with $TV$ yields the result.
\end{proof}

A classical result of toric geometry states that $\TV^\C$ induces an equivalence of the category of fans and the category of normal varieties with actions of $(\C^*)^n$ with dense open orbit and equivariant morphisms. We do not know if the analogous result holds in the spectral setting. 

\begin{Question}
    Is the functor $\TV:\mathrm{Fans}\to \calA ct(SDM)$ fully faithful? Is there an easy recognition principle which elements in $\calA ct(SDM)$ lie in the essential image of $\TV$?
\end{Question}

\begin{remark}
    The action diagram of Remark \ref{actiondiag} is not enough to reconstruct the toric variety, because it already forgets that $\G_m$ is a strict abelian group object.
\end{remark}

In a later section we want to define toric stacks. For that we want to talk about stacky quotients of actions of subgroups of algebraic tori on toric varieties. The next lemma show that this is independent of which definition of action we choose.

\begin{lemma}\label{cofinal}
    Let $\calA ct\vert_{b=1}$ be the full subcategory on objects $(L,b)$ with $b=1$. The functor $\Delta^{op}\to \calA ct\vert_{b=1}$ described above is cofinal\footnote{Sometimes also called right cofinal.}. 
\end{lemma}

\begin{proof}
    We will denote $\calA ct\vert_{b=1}$ by $\calC$. We want to use (the $\infty$-categorical version of) Quillen's Theorem A as in \cite[Theorem 4.1.3.1]{HTT}. In particular, we need to show that for every $d\in \calA ct$ the category $\Delta^{op} \times_{\calC} \calC_{(L,b)/}$ is contractible. First we want to show that the category is connected. \\
    Assume $d=(\Z,1)$. The diagrams 
    
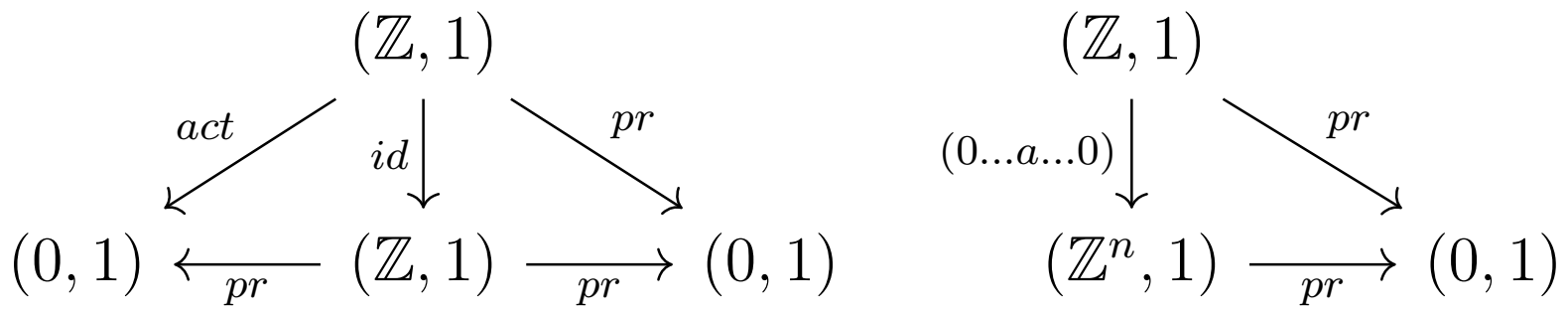
\[\begin{tikzcd}
	& {(\Z,1)} && {(\Z,1)} \\
	{(0,1)} & {(\Z,1)} & {(0,1)} & {(\Z^n,1)} & {(0,1)}
	\arrow["act"', from=1-2, to=2-1]
	\arrow["id"', from=1-2, to=2-2]
	\arrow["pr", from=1-2, to=2-3]
	\arrow["{(0...a...0)}"', from=1-4, to=2-4]
	\arrow["pr", from=1-4, to=2-5]
	\arrow["pr", from=2-2, to=2-1]
	\arrow["pr"', from=2-2, to=2-3]
	\arrow["pr"', from=2-4, to=2-5]
\end{tikzcd}\] show that the action morphism and inclusions of sublattices are both connected to projection in $\Delta^{op} \times_{\calA ct} \calA ct_{(\Z,1)/}$. Moreover, every morphism can be decomposed in actions, inclusions and projections. Hence the category is connected. An analogous analysis shows this also in the case of $d=(\Z^n,1)$, and by equivalence also for all non-trivial lattices $L.$ \\
Fo the case $d=(0,1)$ we can show analogously that everything is connected to the identity.
Next we want to compute the automorphisms of $(\Z,1)\to (0,1)$ defining the projection. 

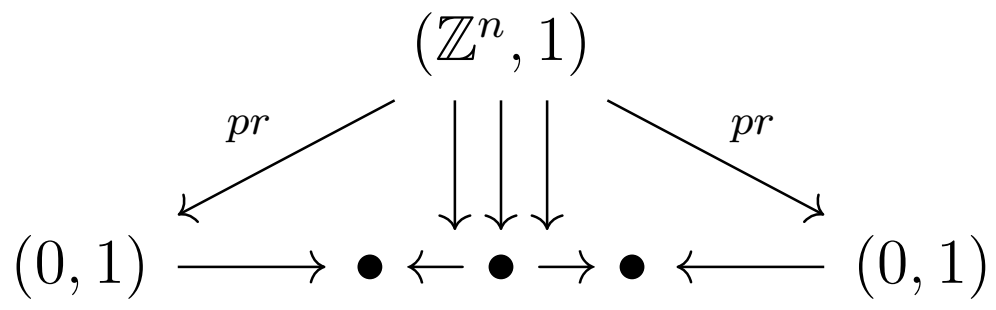
\[\begin{tikzcd}
	& {(\Z^n,1)} \\
	{(0,1)} & {\bullet\leftarrow \bullet\rightarrow \bullet} & {(0,1)}
	\arrow["pr"', from=1-2, to=2-1]
	\arrow[from=1-2, to=2-2]
	\arrow[shift right=3, from=1-2, to=2-2]
	\arrow[shift left=3, from=1-2, to=2-2]
	\arrow["pr", from=1-2, to=2-3]
	\arrow[from=2-1, to=2-2]
	\arrow[from=2-3, to=2-2]
\end{tikzcd}\]

The elements appearing in the bottom line are all of the form $(\Z^n,0)$ and morphisms are induced from $\Delta^{op}$. Since all lattices in the image of $\Delta^{op}$ are canonically in the form $\Z^n$ we can express all morphisms in matrix form. The maps in the image of $\Delta^{op}$ are matrices with all positive entries, since the "inversion" map is not in $\Delta^{op}$. Moreover, we can can factor everything in matrices only containing entries in ${0,1}$, encoding projection, action and inclusion of the identity. After geometric realization the degeneracy maps become homotopy inverses to action and projection. In particular, all automorphisms are homotopic to the identity. \\
This concludes the argument, because $\Delta^{op}$ and $\calA ct\vert_{b=1}$ are one categories, so higher homotopy groups can not appear. 
\end{proof}

\subsection{A note on line bundles}

Classically, the theory of line bundles on toric varieties is very well understood, and especially Cartier divisors compatible with the torus action. 
Later we will need to produce an ample line bundle on a toric variety so we want to quickly recall how to construct those, and how to adapt to the spectral setting.

\begin{definition}[\cite{Wirthmuller}]\label{Chartrip}
    A \emph{characteristic triple} $(\pi,\delta,\lambda)$ for a fan $F$ in a lattice $L$ consits of: \\
    \begin{enumerate}
        \item an epimorphism $\pi:L'\to L$ with rank 1 kernel,
        \item a generator $\delta\in \ker \pi$,
        \item a section $\lambda: F_\to L'$ of $\pi$, such that $\lambda_{\vert \sigma}$ is additive for each $\sigma$.  
    \end{enumerate}
\end{definition}

\begin{remark}
    Usually, line bundles on toric varieties are classified by certain convex functions on polyhedra or fans. 
    The map $\lambda$ should be regarded as the restriction of a graph of such functions. 
    The element $\delta$ is a bookkeeping tool to make this translation canonical, see \cite[Subsection 1.5]{Wirthmuller}. 
    We will follow this approach as it generalizes more readily to the spectral setting. 
\end{remark}

\begin{lemma}
    Denote by $F'_{\lambda}$ the fan in $L'$ generated by $\lambda(\sigma)$.
    The map $\TV_{F'_\lambda}\to \TV_F$ induced by $\pi$ is a $\G_m$-bundle.
\end{lemma}

\begin{proof}
    The map $\TV_\pi$ is weakly semistable, by construction. So by \cref{fiberprods}, we can compute the fiber in {Fans}. 
    For every cone $\sigma\in F$, the fiber is equivalent to $\lambda(\sigma)\times \Z\delta$. Therefore, $\TV_\pi$ is a bundle, with fibers $\G_m$, which is trivialized on each $\TV_\sigma$. 
    By the same argument trivializations are compatible.
\end{proof}

\begin{example}
    Let $F$ be the fan spanned by $(1,0), (0,1), (1,1)$, so that $\TV_F\simeq \P^2$. 
    Define the characteristic triple $(\pi,\delta,\lambda)$ as projection along the third coordinate, $\delta=(0,0,1)$ and $\lambda(1,0)=(1,0,0),\lambda(0,1)=(0,1,0),\lambda(-1,-1,-1)$.
    The line bundle defined by the characteristic triple is $O(1)$.
\end{example}

\begin{prop}
    Let $(\pi,\delta,\lambda)$ be a characteristic triple for the fan $F$ in the lattice $L$, such that for all $\sigma$ there is a $l\in L^\vee$, with $l(\delta)=1$, $l(\lambda)\geq0$ and $l(\lambda(v))=0$ if and only if $v\in\sigma$. 
    The line bundle associated to the characteristic triple over $\TV_F$ relative to $\Z$ is ample. 
\end{prop}

\begin{proof}
    This is \cite[Theorem 1.16]{Wirthmuller} for $\Gamma=0$.
\end{proof}

\section{Constructing the Tate curve}\label{Tate}

The origins of spectral algebraic geometry are inherently connected to generalized elliptic curves: The first spectral stack was $(\overline{\calM}_{ell},\calO^{top})$ constructed by Hopkins-Miller and Goerss-Hopkins. 
In an attempt to better understand this construction, Lurie defined spectral algebraic geometry and reinterpreted (the open substack) $(\calM,\calO^{top})$ as the stack $\Mell^{or}$ the stack of oriented elliptic curves of \cite{Ell2}. In this way, Lurie also lifted $\calE$, the universal elliptic curve, to the spectral setting. 
Still, a part is missing, namely the complement of $\calM_{ell}$ in $\overline{\calM}_{ell}$. 
This part is best understood from a deformation-theoretic point of view. 
Near the compactification point, the universal curve is given as a degeneration of elliptic curves into a nodal curve. The universal example is known as the Tate curve. \\
In this section, we want to construct a spectral algebraic lifting of the Tate curve $T^\heartsuit$ over $\Spec \Z \ppb$.
The Tate curve (or more accurately, an affine open in the curve) is best known as the planar curve cut out by the equation $$y^2+xy=x^3+a_4(q) x+a_6 (q)$$ for explicit power series $a_4,a_6$ in $q$. In spectral algebraic geometry, it is not possible to define schemes by equations. \\
Another way to think about the Tate curve (after inverting $q$) is as a quotient $\G_m/q^\Z$. This is how we will construct the Tate curve. As it turns out, toric geometry gives a nice way to produce this quotient, even spectrally, and its degeneration to $q=0$.\\
As suggested in \cite{Survey}, we will split this construction into a three-step process:\\
\begin{enumerate}
    \item We will be constructing the so-called unfolded Tate curve as a toric variety.
    \item The second step will be to complete the unfolded Tate curve. 
    \item Lastly, we will fold the curve via a certain quotient by group action. 
    \item Algebraize the Tate curve.
\end{enumerate}
This construction essentially agrees with the construction of \cite{DR} and is well known to experts; we nonetheless decided to give it in detail, because we will later use the specific construction to obtain more structure on the Tate-curve.

In the following, $\TV$ is always understood to be based on $\S$, i.e, all toric varieties are (flat) stacks over $\Spec \S$.\\

\emph{Step 1:}For the construction of the Tate curve, we fix the following fan $F_{Tate}$ consisting of the zero cone $\{0\}$  and the families of cones 
\begin{align}
    \sigma_n&=\{ (a,b)\in \Z^2 \mid na\leq b\leq (n+1)a\},\\ 
    \tau_n&=\{(a,b)\in\Z^2 \mid na=b\}
\end{align}  for each integer $n$. We will see that a certain completion of $\TV_{F_{Tate}}$ is the unfolded Tate curve. \\ 

\begin{figure}[h]
    \centering
    \begin{tikzpicture}[scale=0.4]

    \fill (0,0) circle (1pt);
    \fill (0,1) circle (1pt);
    \fill (0,2) circle (5pt);
    \fill (0,3) circle (1pt);
    \fill (0,4) circle (1pt);
    \fill (0,5) circle (1pt);
    \fill (0,6) circle (1pt);
    \fill (0,7) circle (1pt);
    \fill (0,8) circle (1pt);
    \fill (0,9) circle (1pt);
    \fill (0,10) circle (1pt);
    \fill (1,0) circle (3pt);
    \fill [violet] (1,1) circle (3pt);
    \fill (1,2) circle (3pt);
    \fill [red] (1,3) circle (3pt);
    \fill [green] (1,4) circle (3pt);
    \fill [blue] (1,5) circle (3pt);
    \fill [pink] (1,6) circle (3pt);
    \fill (1,7) circle (3pt);
    \fill (1,8) circle (3pt);
    \fill (1,9) circle (3pt);
    \fill (1,10) circle (3pt);
    \fill [violet] (2,0) circle (3pt);
    \fill [violet] (2,1) circle (3pt);
    \fill (2,2) circle (3pt);
    \fill (2,3) circle (3pt);
    \fill [red] (2,4) circle (3pt);
    \fill [red] (2,5) circle (3pt);
    \fill [green] (2,6) circle (3pt);
    \fill [green] (2,7) circle (3pt);
    \fill [blue] (2,8) circle (3pt);
    \fill [blue] (2,9) circle (3pt);
    \fill [pink] (2,10) circle (3pt);
    \fill [violet] (3,0) circle (3pt);
    \fill [violet] (3,1) circle (3pt);
    \fill (3,2) circle (3pt);
    \fill (3,3) circle (3pt);
    \fill (3,4) circle (3pt);
    \fill [red] (3,5) circle (3pt);
    \fill [red] (3,6) circle (3pt);
    \fill [red] (3,7) circle (3pt);
    \fill [green] (3,8) circle (3pt);
    \fill [green] (3,9) circle (3pt);
    \fill [green] (3,10) circle (3pt);
    \fill [violet] (4,0) circle (3pt);
    \fill [violet] (4,1) circle (3pt);
    \fill (4,2) circle (3pt);
    \fill (4,3) circle (3pt);
    \fill (4,4) circle (3pt);
    \fill (4,5) circle (3pt);
    \fill [red] (4,6) circle (3pt);
    \fill [red] (4,7) circle (3pt);
    \fill [red] (4,8) circle (3pt);
    \fill [red] (4,9) circle (3pt);
    \fill [green] (4,10) circle (3pt);
    \fill [violet] (5,0) circle (3pt);
    \fill [violet] (5,1) circle (3pt);
    \fill (5,2) circle (3pt);
    \fill (5,3) circle (3pt);
    \fill (5,4) circle (3pt);
    \fill (5,5) circle (3pt);
    \fill (5,6) circle (3pt);
    \fill [red] (5,7) circle (3pt);
    \fill [red] (5,8) circle (3pt);
    \fill [red] (5,9) circle (3pt);
    \fill [red] (5,10) circle (3pt);
    \fill [violet] (6,0) circle (3pt);
    \fill [violet] (6,1) circle (3pt);
    \fill (6,2) circle (3pt);
    \fill (6,3) circle (3pt);
    \fill (6,4) circle (3pt);
    \fill (6,5) circle (3pt);
    \fill (6,6) circle (3pt);
    \fill (6,7) circle (3pt);
    \fill [red] (6,8) circle (3pt);
    \fill [red] (6,9) circle (3pt);
    \fill [red] (6,10) circle (3pt);
    \fill [violet] (7,0) circle (3pt);
    \fill [violet] (7,1) circle (3pt);
    \fill (7,2) circle (3pt);
    \fill (7,3) circle (3pt);
    \fill (7,4) circle (3pt);
    \fill (7,5) circle (3pt);
    \fill (7,6) circle (3pt);
    \fill (7,7) circle (3pt);
    \fill (7,8) circle (3pt);
    \fill [red] (7,9) circle (3pt);
    \fill [red] (7,10) circle (3pt);
    \fill [violet] (8,0) circle (3pt);
    \fill [violet] (8,1) circle (3pt);
    \fill (8,2) circle (3pt);
    \fill (8,3) circle (3pt);
    \fill (8,4) circle (3pt);
    \fill (8,5) circle (3pt);
    \fill (8,6) circle (3pt);
    \fill (8,7) circle (3pt);
    \fill (8,8) circle (3pt);
    \fill (8,9) circle (3pt);
    \fill [red] (8,10) circle (3pt);

    \draw [thick] (0,0) to (0,10);
    \draw [thick] (-1,2) to (8,2);

    \filldraw [black!80, opacity=0.1] (0,2)  -- (8,10)  -- (8,2)  -- cycle;
    \filldraw [red!80, opacity=0.1] (0,2)  -- (8,10)  -- (4,10)  -- cycle;
    \filldraw [green!80, opacity=0.1] (0,2)  -- (4,10)  -- (2.7,10)  -- cycle;
    \filldraw [violet!80, opacity=0.1] (0,2)  -- (2,0)  -- (8,0)  --  (8,2)  -- cycle;

\end{tikzpicture}
    \caption{A part of the fan $F_{Tate}$}
    \label{FTate}
\end{figure}
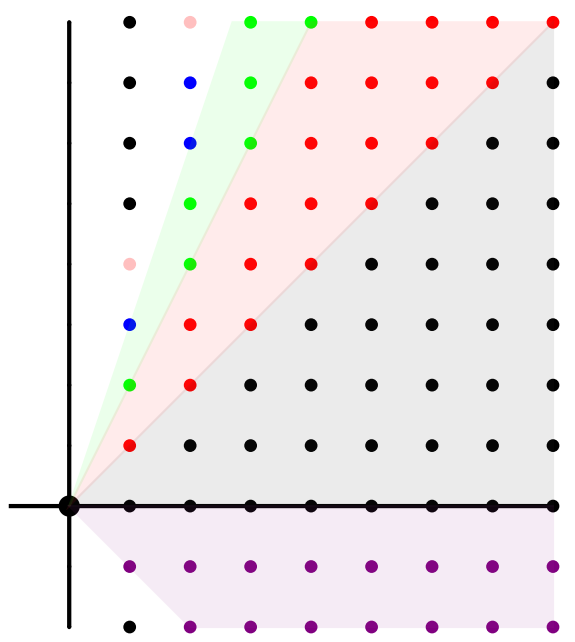

To compute $\TV_{F_{Tate}}$ let us look at the fan morphism $F_{Tate}\to F_0=\{0,\Z_{\geq 0}\}$, given by projection onto the first factor. This morphism gives rise to a map $f:\TV_{F_{Tate}}\to \TV_{F_0}\simeq\A^1$.
First observe that the morphism is weakly semistable, because every cone $\tau_n,\sigma_n$ surjects onto $\Z_{\geq0}$.\\

\begin{lemma}
    The fiber $f^{-1}(\G_m)$ is equivalent to $\G_m^2$. 
\end{lemma}

\begin{proof}
    As we noted the morphism of lattices is weakly semistable. Since $\S$ is noetherian, we can use \Cref{fiberprods}. The claim reduces thus to computing the fiber product $\{0\}\times_{F_0}F_{Tate}$, where the map $\{0\}\to F$ is the inclusion of the 0 cone. 
    The fiber product is therefore given by $0\times \sigma$, for all $\sigma\in F_{Tate}$ mapping to 0, which is just the 0 cone in $F_{Tate}$ itself. 
\end{proof}

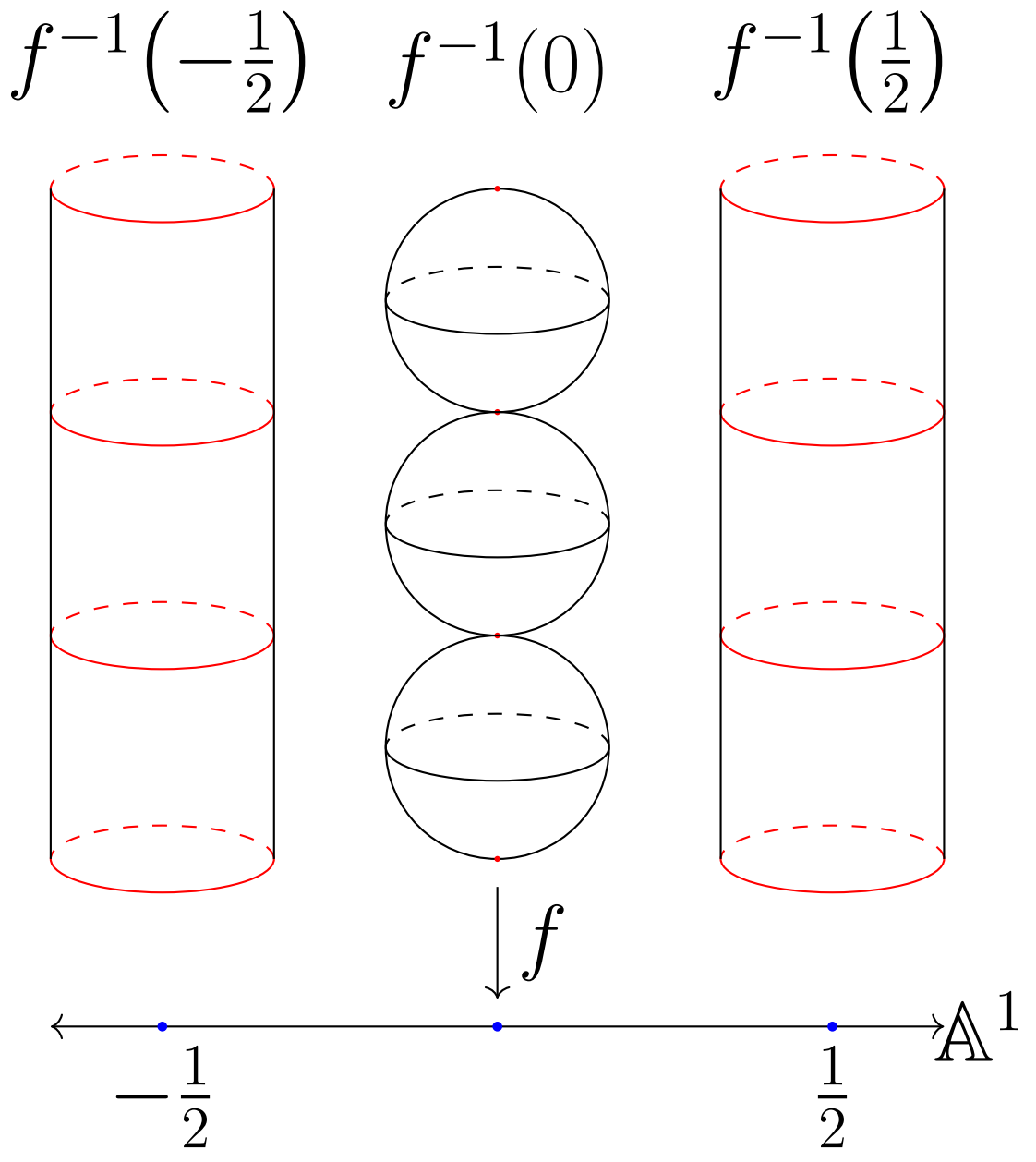
\begin{figure}[h]
    \centering
    \begin{tikzpicture}[scale=0.4]
  \draw (0,0) circle (2cm);
  \draw (-2,0) arc (180:360:2 and 0.6);
  \draw[dashed] (2,0) arc (0:180:2 and 0.6);
  \fill[fill=red] (0,2) circle (1.5pt);
  \draw (0,4) circle (2cm);
  \draw (-2,4) arc (180:360:2 and 0.6);
  \draw[dashed] (2,4) arc (0:180:2 and 0.6);
  \fill[fill=red] (0,-2) circle (1.5pt);
  \draw (0,-4) circle (2cm);
  \draw (-2,-4) arc (180:360:2 and 0.6);
  \draw[dashed] (2,-4) arc (0:180:2 and 0.6);
  \fill[fill=red] (0,-6) circle (1.5pt);
  \fill[fill=red] (0,6) circle (1.5pt);
  \node[above] at (0,7) {\Large$f^{-1}(0)$};

  \draw[color = red] (4,2) arc (180:360:2 and 0.6);
  \draw[dashed, color = red] (8,2) arc (0:180:2 and 0.6);
  \draw[color = red] (4,-2) arc (180:360:2 and 0.6);
  \draw[dashed, color = red] (8,-2) arc (0:180:2 and 0.6);
  \draw[color = red] (4,-6) arc (180:360:2 and 0.6);
  \draw[dashed, color = red] (8,-6) arc (0:180:2 and 0.6);
  \draw[color = red] (4,6) arc (180:360:2 and 0.6);
  \draw[dashed, color = red] (8,6) arc (0:180:2 and 0.6);
  \draw (4,-6) -- (4,6);
  \draw (8,-6) -- (8,6);
  \node[above] at (6,7) {\Large$f^{-1}\bigl(\frac{1}{2}\bigr)$};
  
  \draw[color = red] (-8,2) arc (180:360:2 and 0.6);
  \draw[dashed, color = red] (-4,2) arc (0:180:2 and 0.6);
  \draw[color = red] (-8,-2) arc (180:360:2 and 0.6);
  \draw[dashed, color = red] (-4,-2) arc (0:180:2 and 0.6);
  \draw[color = red] (-8,-6) arc (180:360:2 and 0.6);
  \draw[dashed, color = red] (-4,-6) arc (0:180:2 and 0.6);
  \draw[color = red] (-8,6) arc (180:360:2 and 0.6);
  \draw[dashed, color = red] (-4,6) arc (0:180:2 and 0.6);
  \draw (-8,-6) -- (-8,6);
  \draw (-4,-6) -- (-4,6);
  \node[above] at (-6,7) {\Large$f^{-1}\bigl(-\frac{1}{2}\bigr)$};

    \draw[<->] (-8,-9) -- (8,-9);
    \fill[fill=blue] (0,-9) circle (2.5pt);
    \fill[fill=blue] (-6,-9) circle (2.5pt);
    \fill[fill=blue] (6,-9) circle (2.5pt);
    \node[below] at (-6,-9) (a) {\Large$-\frac{1}{2}$};
    \node[below] at (6,-9) (b) {\Large$\frac{1}{2}$};
    \node at (8.5,-9) (c) {\hspace{3pt}\Large$\A^1$};
    \draw[->] (0,-6.5) to node[right] {\Large$f$} (0,-8.5);
    
\end{tikzpicture}
    \caption{Depiction of $\TV_{F_{Tate}}$ over the complex numbers}
    \label{unfolded}
\end{figure}

The more interesting case is the fiber over $0\in \A^1$. A depiction can be found in \Cref{unfolded}. 
The more geometrically minded reader might like to compare the following computations to the figure. \\
By definition, for every cone $\sigma_n\in F_{Tate}$, the toric variety $\TV_{\sigma_n}$ is an affine 2-space $\A^2.$ 
All of them are glued along $\TV_{\tau_n}\simeq \A^1\times \G_m$. 
Since in every copy of $\A^2$ the preimage of $\G_m$ is $\G_m^2$, we might expect the coordinate axes to get mapped to $0$. This is indeed true in the following sense:
Fix $\sigma_n\in F_{Tate}$. The dual cone $\check{\sigma}_n$ is generated by $\{(-n,1),(n+1,-1)\}$, where we identify $(\Z^2)^\vee\cong \Z^2$ using the standard bilinear form on $\Z^2$. 
We obtain coordinates $$\A^2\simeq \Spec \S[x^{-n}y ,x^{n+1}y^{-1}] \xrightarrow{f}\A^1\simeq\Spec \S[x].$$ In these coordinates $f$ is given by the map $$\S[x]\to \S[x^{-n}y,x^{n+1}y^{-1}], x\mapsto x=x^{-n}y\cdot x^{n+1}y^{-1}.$$ 
The fiber we are looking for can thus be described by $\S\otimes_{\S[x]} \S[x^{-n}y,x^{n+1}y^{-1}]$, with $x$ acting as 0 on the left factor. 
On $\pi_0$ this is just $$\Z\otimes_{\Z[x]}\Z[x^{-n}y,x^{n+1}y^{-1}]\cong \Z[q,w]/(qw).$$ Therefore, the fiber is a lift of the coordinate axes to $\S$.\\

Let us now turn our attention to $\tau_n$, generated by $(1,n)$. The dual $\check\tau_n$ is generated by $$\{(-n,1),(n+1,-1),(n,-1)\}$$ or alternatively $$\{\{(-n,1),(-n+1,1),(n,-1)\}.$$ 
Therefore we obtain coordinates $$\A^1\times \G_m\simeq \Spec \S[x^{-n}y,x^{n+1}y^{-1},x^ny^{-1}]\simeq \Spec \S[x^{-n}y,x^{-n+1}y,x^ny^{-1}].$$ 
The inclusion of $\tau_n$ as bottom face of $\sigma_n$ is then given as localization $$\S[x^{-n}y,x^{n+1}y^{-1}]\to\S[x^{-n}y,x^{n+1}y^{-1},x^ny^{-1}],$$ while the inclusion of $\tau_n$ as top face of $\sigma_{n-1}$ induces $$\Spec \S[x^{-n+1}y,x^ny^{-1}]\to\Spec \S[x^{-n}y,x^{-n+1}y,x^ny^{-1}].$$ 
Recall that we need to consider the pushout $\TV_{\sigma_n}\xleftarrow{} \TV_{\tau_n}\to \TV_{\sigma_{n-1}}$.
The generic fiber, meaning over $\G_m\subseteq \A^1$, is just $\G_m^2$. 
The special fiber is given by two sets of coordinate axis pairs identified along a $\G_m$, where the $\G_m$ is included in one axis in the canonical way.
More generally, in both cases, the inclusions factor through one axis, and if we restrict our pushout to these components, we recover the construction of $\P^1$.  
The whole stack $\TV_{F_{Tate}}$ is thus given by a $\Z$ indexed set of coordinate systems, say $C_n$ where one axis $\A^1\subseteq C_n$ is glued to one axis of $C_{n-1}$, while the other copy $\A^1\subseteq C_n$ is glued to $C_{n+1}$. 
This shows that the fiber over $0\in \A^1$ is given as $\Z$ indexed family of $\P^1$, where the pole $\infty$ in the $n$-th copy of $\P^1$ is glued to the pole $0$ in the $n+1$st copy of $\P^1$.

\Cref{unfolded} shows a picture of $\TV_{F_{Tate}}$ over the complex numbers. 
Each sphere corresponds to a cone $\tau_n$, and each pinching point (drawn in red) corresponds to a cone $\sigma_n$. In particular, the variety $\TV_{\sigma_n}$ is the variety where you only pinch one of the red circles in $\G_m$ to a point (the resulting intersection are the coordinate axes described above). \\

\emph{Step 2:} We define the stack $\widehat{\TV}_{F_{Tate}}$ as the completion of $\TV_{F_{Tate}}$ along the fiber $f^{-1}(0)$, over $\Spf \S\llb q\rrb$. \\

\emph{Step 3:} Consider the morphism $\nu:\Z^2\to \Z^2, (a,b)\mapsto (a,a+b)$. This morphism defines a group action of $\Z$ on the fan $F_{Tate}$ so $\Z$ acts via $\TV\nu$ on $\TV_{F_{Tate}}$ and $\widehat{\TV}_{F_{Tate}}$.
Moreover, the group acts on $F_{Tate}$ over $F_0$ so we obtain an action on $\TV_{F_{Tate}}$ as a stack over $\A^1$.\\ 
The action on $F_{Tate}$ shifts the cones $\tau_n\mapsto \tau_{n+1}$ and $\sigma_n\mapsto \sigma_{n+1}$. 
Therefore, the action on $\TV_{F_{Tate}}$ is easy to describe in the special fiber. It shifts the copies of $\P^1$ we identified in Step 1.\\
On the generic fiber, the action is given by $\S[q^\pm,w^\pm]\to \S[q^\pm,w^\pm], q\mapsto q, w\mapsto qw$, where $q$ is the coordinate of the target of $\TV_{F_{Tate}}\to \A^1$. 
In particular, this map is the shearing map of $\G_m$ acting on itself via left multiplication.
In particular, over a $(R-)$ point $q$ in $\G_m$ the action of $\nu$ is given by multiplication by $q$.\\ 
This $\Z$ action on $\TV_{F_{Tate}}$ is not free: $q\in \G_m(R)$ can have finite order. 
The action becomes free, however, after base change to $\widehat{\TV}_{F_{Tate}}$.\\

\begin{definition}
    We define the formal scheme $\widehat{T}:=(\widehat{\TV}_{F_{Tate}})/\Z$ (where we take the quotient in formal spectral schemes) over $\Spf \S\bbq$ and call $T$ the \emph{formal Tate curve}. 
\end{definition}

\begin{prop}
    The formal scheme quotient $\widehat{T}$ exists. 
\end{prop}

\begin{proof}
    We first construct the quotient $(\widehat{\TV}_{F_{Tate}})/2\Z$ in formal spectral Deligne-Mumford stacks which exists, because of \cite[Proposition 8.7.6]{Rezksag}.
    Let us consider the subfan of $F_{Tate}$ generated by $F'=\{(\sigma_n)_{-2<n<2}\}$. The $q$-completed toric variety associated to $F'$ has only finitely many proper irreducible components. 
    Each of which include into $(\widehat{\TV}_{F_{Tate}})/2\Z$ as a Zariski open, and together they cover the quotient.
    In particular, the quotient is a scheme. 
    The statement follows as $\widehat{T}\simeq ((\widehat{\TV}_{F_{Tate}})/2\Z)/(\Z/2/Z)$.
\end{proof}

\textit{Step 4:} Define the characteristic triple  $(\pi,\delta,\lambda)$ as $\pi:\Z^3\to \Z^2$ projection along $(0,0,1)$, $\delta=(0,0,1)$, $\lambda(1,i)=(1,-i,-i(i+1)/2)$.
As in subsection 2.4, this defines a line bundle on $\widehat{\TV}_{F_{Tate}}$. 
In section 8.4 of \cite{Dirichlet} the authors compute the $(\Z$-invariant) sections of the line bundle on $(\widehat{\TV}_{F_{Tate}}^\heartsuit)$, called $\theta$-functions, and show
that this line bundle is relatively ample. 

\begin{definition}
    Let $\mathfrak{X}$ be a (formal) spectral scheme over a (formal) base scheme $\Spec A$. 
    We say a line bundle $\calL$ over $X$ is relatively ample, if $\pi_0\calL$ is relatively ample over $X^\heartsuit$ relative to $S^\heartsuit$.
\end{definition}

\begin{thm}[Langwen Hui]
    Let $\mathfrak{X}\to \Spf A$ be a proper formal scheme over some connective, noetherian complete $\E_{\infty}$-ring $A$. 
    Let $\calL$ be a relatively ample line bundle over $\mathfrak{X}\to \Spf A$. Then $\mathfrak{X}$ algebraizes to a stack $X\to \Spec A$.
\end{thm}

\begin{proof}
    This proof will appear in \cite{Langwen}.
    Let us give a sketch of the proof.\\ 
    By \cite[8.5.0.3]{SAG} it suffices to produce a closed immersion $\mathfrak{X}\into \widehat{\P}^N$ into some \emph{smooth} large projective space over $\Spf A$.
    We do this by showing that there is a localy split map $\calL^{\otimes n}\to \calO_{\mathfrak{X}^{N+1}}$.
    Because $\calL$ is assumed to be ample, we know that $\pi_0\calL=:L$ is ample over $\mathfrak{X}^\heartsuit$. 
    In particular, for some $m\gg 0$, there are sections \[ s_0,...,s_N\in H^0(\mathfrak{X}^\heartsuit, L^{\otimes m})\] defining a closed immersion $\mathfrak{X}^\heartsuit\into (\P^N)^\heartsuit$. 
    Analyzing the descent spectral sequence for $\Gamma(\mathfrak{X},\calL^m)$, we can see that for large $m$, the map $\pi_0\Gamma(\mathfrak{X},\calL^{\otimes m})\to H^0(\mathfrak{X}^\heartsuit, L^{\otimes m})$ becomes an isomorphism. 
    This provides lifts of $s_0,...,s_N: \calO^{N+1}\to \calL^{\otimes m}$, whose dual we were looking for.   
\end{proof}

\begin{cor}
    The formal Tate curve $\widehat{T}:=\widehat{\TV}_{F_{Tate}}/\Z \to \Spf \S\llb q\rrb$ algebraizes to a curve $T\to \Spec \S \llb q \rrb$. Furthermore, $T^\heartsuit$ is the classical Tate curve.
\end{cor}

\begin{proof}
    We want to construct an ample line bundle on $\widehat{T}$.
    Consider the characteristic triple $(\pi,\delta,\lambda)$ as $\pi:\Z^3\to \Z^2$ projection along $(0,0,1)$, $\delta=(0,0,1)$, $\lambda(1,i)=(1,-i,-i(i+1)/2)$ on $F_{Tate}$. 
    As in section 2.4, this defines a line bundle $\calL$ on $\TV_{F_{Tate}}$.
    In \cite[Section 8.4]{Dirichlet}, the authors compute the $\Z$-invariant sections on $\widehat{\TV}_{F_{Tate}}$ over $\C$, and show that this line bundle descends to an ample line bundle on $\widehat{\TV}_{F_{Tate}}/\Z$.
    A less detailed, but integral analogue can be found in \cite[VII.1.10]{DR}. 
    The integral line bundle is precisely $\pi_0\calL$ on  $(\widehat{\TV}_{F_{Tate}}/\Z)^\heartsuit$. 
    Therefore, by above theorem, we know that $\widehat{T}$ algebraizes over $\Spec \S \llb q\rrb$. 
    Furthermore, we can do the same procedure based on $\Z$, instead of $\S$. 
    Flatness of toric varieties guarantees that basechanging along $\S\to \Z$ agrees with $(-)^\heartsuit$. 
    Furthermore, because formal completion and taking the $\Z$-quotient are colimits they commute with colimits. 
    Because algebraization is compatible with base change we get the claim.\\
    Alternatively, one can obsereve that the above defined line bundle is strict, and one can define a spectral Proj construction. This takes more effort and will appaer in \cite{Langwen}.
\end{proof}

Later, we will also examine the smooth locus of $T$. It has a similar description.
Consider the fan $\bar{F}$ consisting of the cones $\tau_n$ and $\{0\}$, i.e. the fan obtained from $F_{Tate}$ by removing the 2-dimensional cones. The same calculation as before shows that $TV_{\bar{F}}$ is given by the unfolded Tate curve without the singularities. \\
The discussion above also shows that we can compute the smooth locus $T^{reg}$ of $T$ by taking a completion and quotient of $TV_{\bar{F}}$.

\begin{thm}
    The universal deformations $\calN$ of the $n$-gons lift over the sphere.
\end{thm}

\begin{proof}
    The construction is analogous to the construction of the Tate curve, but replacing the quotient $(\widehat{\TV}_{F_{Tate}})/\Z$ by $(\widehat{\TV}_{F_{Tate}})/n\Z$. 
\end{proof}

\section{Group Structure on the Tate curve}

We begin this section by reminding the reader of the structure of the classical Tate curve $T^\heartsuit$. The Tate curve is the universal deformation of the nodal curve, with generic fibers given by elliptic curves. In particular, the generic fibers carry group structures. This group structure extends to the smooth locus $(T^\heartsuit)^{reg}$.
Furthermore, the group that is given by the smooth locus of the special fiber acts on the entire nodal curve.
In spectral algebraic geometry, one has to differentiate between elliptic curves and strict elliptic curves. For elliptic cohomology, the elliptic curves must be strict (see \cite[Chapter 7]{Ell2}. 

\subsection{Towards the group structure on the Tate curve}
In this section, we show that the spectral Tate curve $T$ constructed above has similar properties over the sphere. 

Consider the following morphism of lattices $$\Z^3\to \Z^2, (a,b,c)\mapsto (a,b+c).$$ 
This can be thought of as an addition map of $\Z^2$ over $\Z$, since $\Z^3\cong \Z^2 \times_\Z \Z^2$. 
Now looking at $\overline{F}=\{ \{0\},\tau_n\}$ in $\Z^2$ this yields a morphism of fans $$\overline{F}\times_{F_0} \overline{F}\to \overline{F}$$ and also $\overline{F}\times_{F_0} F\to F$.

As can be seen in \cite[Lemma 2.2.6]{Pullback} taking this pullback commutes with applying the functor $TV_-$. 
So the smooth locus of $TV_F$ obtains the structure of an abelian group object and it acts on $TV_F$. \\
Our interest though lies in the quotient $T=\widehat{TV}_F/\nu$ over the disc $\Spec\S\llb q\rrb$. The next lemma describes $T$ over the punctured disc.

\begin{lemma}
    The smooth locus of the unfolded Tate curve is a group scheme, with formal group canonically equivalent to $\widehat{\G}_m$. Moreover, the group structure descends to $\widehat{T^{reg}}/\Z$.
\end{lemma}

\begin{proof}
    Consider the map of fans $\{0\} \times \Z_{\geq 0} \to \overline{F}$ over $\Z_{\geq 0}$. 
    This map includes $\G_m$ into $\TV_{\overline{F}}$. 
    Moreover, this inclusion is compatible with the group structure obtaines by $\Z^3\to \Z^2, (a,b,c)\mapsto (a,b+c)$ and therefore induces a map of group stacks. 
    The calculation of the previous section shows that whenever $q$ is invertible, this inclusion is an equivalence and when $q=0$, this includes $\G_m$ as the identity component of $\G_m\times \Z$.
    In particular, this inclusion induces an equivalence after completion at the unit.\\
    Furthermore, we can see that whenever $q$ is nilpotent, the $\Z$ action is equivalent to the translation of $\G_m\times \Z$ in the $\Z$ factor. 
    Therefore, after completion of $\TV_{\overline{F}}$, the quotient identifies the $\G_m$ copies and the group structure descends. 
\end{proof}

\begin{warning}
    We do not know if this group structure algebraizes! 
    The smooth locus is no longer proper, so the Grothendieck existence theorem no longer applies. \\
    Classiclally, this can be circumvented as a (possibly degenrate) elliptic curve with a section admits a unquie group structure on the smooth locus acting on the whole curve. 
    Spectrally, this is no longer true.
\end{warning}

\begin{Question}
    Does the group structure on $(T^\heartsuit)^{reg}$ lift to a strict group structure on $T^{reg}$? 
\end{Question}

The above computations show that if this group structure lift as the algebraization of group structure on the formal scheme, then the formal group has to be equivalent to $\widehat{\G}_m$.

Such a group structure would give rise to a definition of equivariant Tate $K$-theory, extending the definition of \cite{SmoothK}.

\subsection{Compactifications of $\Mell$}

We assume for the remainder of this subsection that the group structure of the formal Tate curve algebraizes.
By the above analysis of the group structure it follows that choosing a preorientation or an orientation for the Tate curve is equivalent to choosing the same for $\widehat{\G}_m$.\\
In \cite{Ell2} it was shown that the universal cases for a preorientation and an orientation of $\G_m$ are given by $\Sig\C P^\infty$ (\cite[Proposition 4.3.25]{Ell2}) and $KU$ (\cite[Corollary 4.3.27, Theorem 6.5.1]{Ell2}). 
So basechanging from $\S$ to these cases we get the universal preoriented and universal oriented Tate curves living over $\Spec \Sig \C P^{\infty}\llb q \rrb $ and $\Spec KU\llb q \rrb$ respectively. \\

Furthermore, the (oriented, preoriented) Tate curve can be used to compactify the respective moduli stacks constructed by Lurie. 
The reflection along the cone $\tau_0$ in the Fan $F_{Tate}$ induces the multiplication by $-1$ map, which we can quoteint by to obtain a curve over $\Spec \S\llb q\rrb/C_2$ and the respective basechanges. 
We obtain the following spans, where the top horizontal arrows classify the constructed Tate curves and the left vertical arrows, are canonical. \\

\begin{center}
\begin{tikzcd}
    \Spec \S\ppb/C_2 \ar{r}\ar{d}& \mathcal{M}_{ell}\\
    \Spec \S\llb q \rrb/C_2 & 
\end{tikzcd}
\end{center}
\begin{center}
\begin{tikzcd}
    \Spec \Sig \C P^{\infty} \ppb/C_2 \ar{r}\ar{d}& \mathcal{M}^{pre}_{ell}\\
    \Spec \Sig \C P^{\infty}\llb q \rrb/C_2 & 
\end{tikzcd}
\begin{tikzcd}
    \Spec KU\ppb/C_2 \ar{r}\ar{d}& \mathcal{M}^{or}_{ell} \\
    \Spec KU\llb q \rrb/C_2 & 
\end{tikzcd}
\end{center}

\begin{warning}\label{technical}
    In \cite{DMC} the authors construct compactifications of the stacks of genus $g$ curves. 
    Similarly, Lurie suggests in \cite[Section 4.3]{Survey} to use the pushout of the latter of the above spans as definitions for the compactification $\Mellb^{or}$ of $\Mell^{or}$. 
    Analogously, one can use the other two spans to define $\overline{\calM}_{ell},\overline{\calM}_{ell}^{pre}$. 
    We do not claim the existence of any of the above pushouts. 
    A prove of the existence of the last diagram appears in \cite{SmoothK} for in 1-localic stacks. 
\end{warning}

\begin{thm}
    The pushouts in 1-localic stacks exist.
\end{thm}

\begin{proof}
    For the oriented case this is \cite[Lemma 6.8]{SmoothK}. 
    For the other cases the argument can be adapted by analyzing the Mayer-Vietoris sequence and using connectivity. 
\end{proof}

Now looking at the curves over this pushout we also obtain universal generalized elliptic curves. For example:

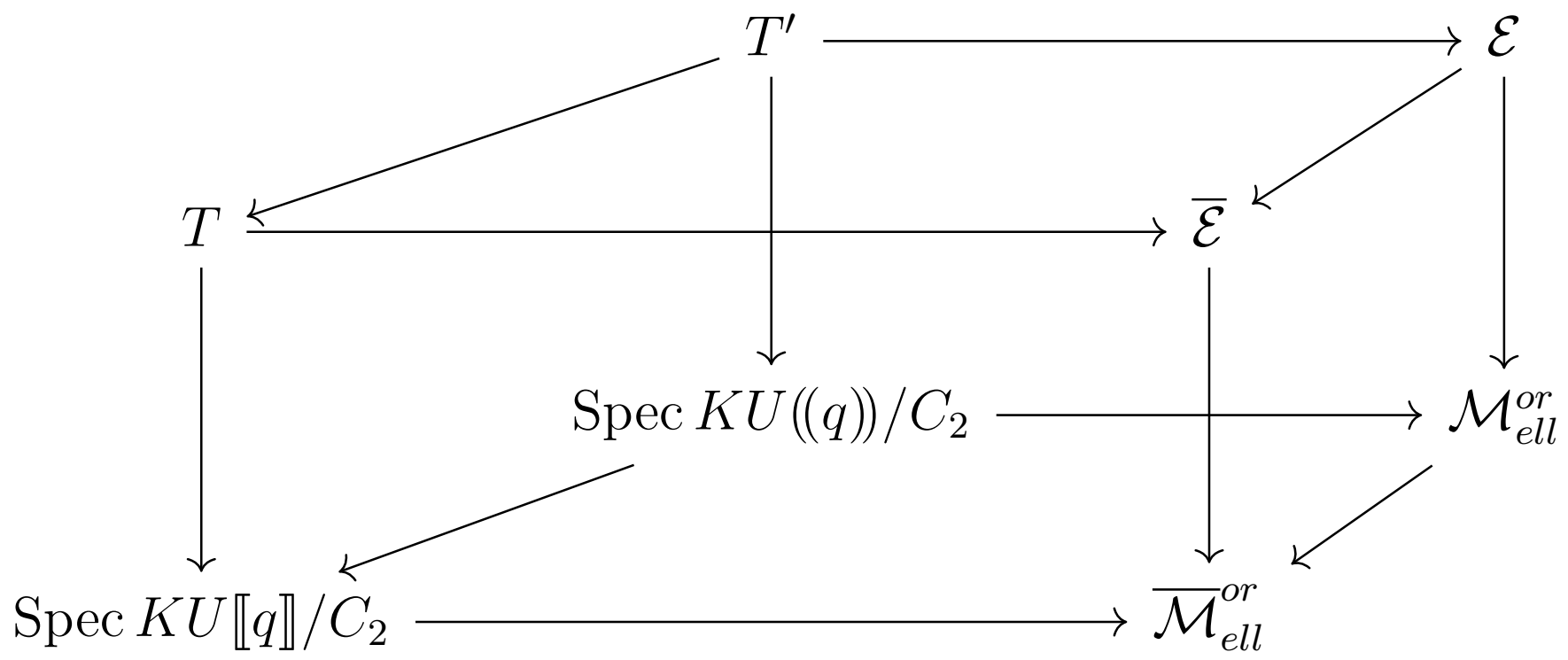
\[\begin{tikzcd}
	& {T'} && \mathcal{E} \\
	T && \overline{\mathcal{E}} \\
	& \Spec KU\ppb/C_2 && \mathcal{M}_{ell}^{or} \\
	\Spec KU \llb q \rrb/C_2 && \overline{\mathcal{M}}_{ell}^{or}
	\arrow[from=1-2, to=1-4]
	\arrow[from=1-2, to=2-1]
	\arrow[from=1-2, to=3-2]
	\arrow[from=1-4, to=2-3]
	\arrow[from=1-4, to=3-4]
	\arrow[from=2-1, to=2-3]
	\arrow[from=2-1, to=4-1]
	\arrow[from=2-3, to=4-3]
	\arrow[from=3-2, to=3-4]
	\arrow[from=3-2, to=4-1]
	\arrow[from=3-4, to=4-3]
	\arrow[from=4-1, to=4-3]
\end{tikzcd}\]

\begin{remark}
	One could try to apply similar techniques to construct (a candidate for) $\calM_{cub}^{or}$ by adding the universal deformation of the cuspidal curve, however as in \cite[Proposition 6.26]{versal} one can see that such a deformation does not exist for the cuspidal curve at primes 2 and 3. 
    More generally, \cite[Proposition 1.6.20]{Ell2} shows that $\calM_{cub}$ should not exist before going to the oriented case or at least should not behave in the expected way (as one can always forget orientations). 
\end{remark}

\begin{remark}
    Using the $n$-gons one can construct similar compactifications of $\Mell$. These have been labeled $\mathcal{E}ll_m$ in \cite{Cesnavicius}. 
    One nice feature is that (after inverting $n$) $\Mellb^\heartsuit$ is etale over $\mathcal{E}ll_m$, so by etale rigidity one has lifts of elliptic curves with full level structures. 
    Such lifts were also constructed in \cite{HL} using log geometry. We expand on this discussion and precise comparisons in future work. 
\end{remark}

\section{Toric stacks}\label{stacks}

We have seen that all geometric objects defined were schemes. However, there are also many interesting stacky examples of geometric objects with torus actions. We are particularly interested in $\A^1/\G_m$ and weighted projective spaces $\P(a_1,..,a_n)$. The latter appears naturally in the theory of equivariant elliptic cohomologies. 

In this section, we want to introduce the theory of toric stacks as in \cite{toricstacks} encompassing the theory of toric varieties, weighted projective spaces, and $\A^1/\G_m$. 

Generally, the groups we want to quotient by are not discrete, and as such, the theory of spectral Deligne-Mumford stacks is insufficient. We therefore want to recall the definition of geometric stacks. 

\begin{definition}
    A \emph{geometric stack} is a $fpqc$ sheaf $X:\CAlg\to \calS$, such that the diagonal $X\to X\times X$ is affine and such that there exists an $\E_\infty$-ring $A$ admitting a faithfully flat map $\Spec A\to X$.
\end{definition}

\begin{remark}
    Geometric stacks have, for example, been studied in \cite{RokMP} under the name of non-connective geometric stacks. 
\end{remark}

\begin{definition}
    A \emph{stacky fan} is a pair $(F,\beta: L\to N)$ consisting of a fan $F$ in the lattice $L$ and a lattice morphism $\beta: L\to N$ with finite cokernel.
\end{definition}

The morphism $\beta$ defines a subgroup of $\G_m^n=\Spec \S[L^\vee]$ as the fiber of $$\Spec(\S[\beta^\vee]):\Spec \S[L^\vee]\to \Spec \S[N^\vee],$$ which we will denote by $G_\beta$. This subgroup acts on $\TV_F$ as in Remark \ref{actiondiag}. 

\begin{definition}
    Let $(F,\beta:L\to N)$ be a stacky fan with associated group $G_\beta.$ We define the toric stack $$\mathfrak{T}_{F,\beta}:=\TV_F/G_\beta:=\colim_{\Delta^{op}} \TV_F\times G_\beta^{\bullet},$$ where we take the colimit in geometric stacks. 
\end{definition}

\begin{remark}
    Recall that by Lemma \ref{cofinal} the different notions of actions we considered lead to equivalent quotients. In particular, the definition of toric stacks does not change if we strictify our action to $\calA ct.$ 
\end{remark}

\begin{numex}
    \begin{enumerate}
        \item Let $F=\Z_{\geq0}\subseteq \Z$ be the fan defining $\A^1$ and let $\beta:\Z\to \{0\}$ be the projection. The subgroup inclusion $G_\beta \subseteq\G_m$ is an equivalence in this case. The toric stack $\mathfrak{T}_{F,\beta}$ is thus given by $\A^1/\G_m$. 
        \item Let $F$ be the fan in $\Z^{n+1}$ consisting of the cones $\langle e_i\rangle$ and $\{0\}$, so that $\TV_F\simeq \A^{n+1}\setminus 0.$ Let furthermore $q_0,...,q_n$ be positive coprime integers. We define $\beta:\Z^{n+1}\to \Z^n,$ by $$\beta(e_0)=\sum_{i=1}^n q_0e_i, \hspace{1cm} \beta(e_i)=-q_i e_i.$$ The dual $\beta^\vee: \Z^n\to \Z^{n+1}$, thus is given by $\beta(e_i)=-q_ie_o+q_0e_i$ and has cokernel $\Z^{n+1}\to\Z, e_i\mapsto q_i$. The toric stack $\mathfrak{T}_{F,\beta}$ therefore is equivalent to $(\A^{n+1}\setminus 0)/ \G_m$, where $\G_m$, acts via the weights $(q_0,...,q_n)$, i.e. it is the weighted projective space.  
    \end{enumerate}
\end{numex}

\begin{numrmk}
    \begin{enumerate}
        \item It is also possible to recover the weighted projective spaces with $(q_0,...,q_n)$ not coprime. This would lead to the theory of non-strict toric stacks in \cite{toricstacks}. 
        \item Alternatively, one could have constructed the weighted projective space as a quotient of the unweighted projective space. This construction is given by the fan $F$ of Example \ref{Projective} and $\beta$ a certain scaling map. The image $\beta(F)\subseteq \Z^n$ then recovers the coarse moduli of the weighted projective space as constructed in \cite[page 35]{Fulton}.
    \end{enumerate}
\end{numrmk}

\section{Arithmetic toric varieties}\label{arithmetic}

Another different kind of flavor of toric varieties are the so called arithmetic toric varieties. 
Arithmetic toric varieties describe toric vatrieties for non-split tori (usually considered over non separably closed fields).  
A general account of the idea used here is given in \cite{Galois} and in \cite{twist} for our applications to toric geometry. 
However, we will not use the bulk of the theories developed there. \\

\begin{definition}
	We call a spectral scheme $Y$ over $R$ a twisted or arithmetic toric varity if for some faithful Galois extension $R\to S$ the pullback $Y\times_{\Spec R} \Spec S$ is a toric variety over $S$. 
\end{definition}

Let $R\to S$ be a faithful finite Galois extension of $\E_\infty$-rings with Galois group $G$ as in \cite{RognesGalois}. Let $X$ be a toric variety over $S$ associated to a fan $F$ in a lattice $L$. 

\begin{prop}
	Let $G$ act on the fan $F$. 
    This action gives rise to an equivariant morphism $\TV_F\to \Spec S$, which exhibits $(\TV_F)/G\to \Spec R$ as an arithmetic toric variety, where we take the quotient $(\TV_F)/G$ in $SpSch$. 
\end{prop}

\begin{proof}
	The action of $G$ on the fan $F$ gets carried to an action on $\TV_F^S$. 
    We define the morphism $\TV_F^S\times_{\Spec S} \Spec S\to \Spec S$ as the projection, where we let $G$ act on $\Spec S$ via the Galois action. 
    The quotient $(\TV_F^S)/G$ exists essentially by \cite{Rydh}. First choose a $G$-stable Zariski covering. This allows us to construct the quotient locally, where it exists as the scheme of fixed points of the global sections.
    In particular, this morphism is equivariant and the quotient $\Spec S/G$ can be canonically identified with $\Spec R$. 
    That this construction indeed yields a twisted form is part of Galois descent established (locally) in \cite[Theorem 9.4]{AkhilGalois}. The global statememt follows again from the compability of covers and pullbacks.
\end{proof}
    
\begin{example}
	Let $F$ be the 0 fan in $\Z$ and let $C_2$ act on $\Z$ by multiplication with $-1$. The fan $F$ defines $\G_m$ over $\Spec \C$ and the resulting arithmetic toric variety over $\R$ is the circle $\Spec \R[x,y]/(x^2+y^2=1)$, while the trivial  $C_2$ action on $\Z$ defines $\R\setminus\{0\}$.
\end{example}

As in the example we see that the 0 fan supports an obvious $C_2$ action. So we can look for Galois extensions of $\E_\infty$-rings realizing this Galois action. As it turns out there is a canonical candidate. 

\begin{example}
	Let $\G_m$ be the universal oriented strict multiplicative group over $KU$ and choose an identification with $\TV_{\{0\}}$ over $\Spec KU$. 
    As in the example above the fan $\{0\}$ gets acted on by a $C_2$. Similarly, $KU$ is a $C_2$ Galois extension of $KO$. The proposition shows that there is a twisted variant of $\G_m$ living over $\Spec KO$, which can be viewed as the universal twisted oriented multiplicative group scheme. 
    Alternatively, one can take the quotiens in SDM to obtain a the twisted $G_m$ over $\Spec KU\llb q\rrb/C_2$, which satisfies better universality properties but is no longer schematic.   
\end{example}

\medskip
\bibliographystyle{amsalpha}
\bibliography{toriclit.bib}

\end{document}